\documentclass{amsart}
\usepackage{mathtools} 
\usepackage{graphics}

\newtheorem{theorem}{Theorem}[section]
\newtheorem{lemma}[theorem]{Lemma}
\newtheorem{corollary}[theorem]{Corollary} 

\theoremstyle{definition}
\newtheorem{assumption}[theorem]{Assumption}
\newtheorem{definition}[theorem]{Definition}

\theoremstyle{remark}
\newtheorem{example}[theorem]{Example}
\newtheorem{remark}[theorem]{Remark}

\newcommand{\mysection}[1]{\section{#1}
      \setcounter{equation}{0}}
\numberwithin{equation}{section}
 
\newcommand\loc{\textnormal{loc}}
\newcommand{\shharp}{=\kern -.5em\|}
\newcommand{\vsharp}{\asymp\kern -.5em\|}

\newcommand{\osc}{\operatornamewithlimits{osc}}

\renewcommand{\eqref}[1]{\text{\rm(\ref{#1})}}

\makeatletter 
 \def\dashint{\operatorname{\,\,\,\mathclap{\!\int}\! \!\text{\bf--}\!\!}}
 \makeatother

\def\dashnorm{\,\,\text{\bf--}\kern-.5em\|}

\newcommand\sfd{{\sf d}}

\newcommand\bB{\mathbb{B}}
\newcommand\bC{\mathbb{C}}

\newcommand\bR{\mathbb{R}}

\newcommand\bS{\mathbb{S}}

\newcommand\cF{\mathcal{F}}

\newcommand\cL{\mathcal{L}}

\newcommand\cO{\mathcal{O}}

\begin{document}

\title[Weak uniqueness
for SDEs with singular coefficients]{Once again about weak uniqueness
for SDE with singular coefficients}

\author[N.V. Krylov]{N.V. Krylov}%
\address{127 Vincent Hall, University of Minnesota,
Minneapolis,
       MN, 55455, USA}
\email{nkrylov@umn.edu}


\begin{abstract} 
We prove  weak uniqueness for
admissible solutions of It\^o's
equations with uniformly nondegenerate $a$ which is almost in VMO
and $b$   in a Morrey class of functions
with low integrability property. If 
$b\in L_{d}$ any solution is admissible.
\end{abstract}

\maketitle

\mysection{Introduction}

Fix a   $\delta\in(0,1]$.
 In this paper  we are dealing with a Borel
function $a$ given on $\bR^{d}=\{x=(x^{1},...,x^{d}),x^{i}\in\bR\}$  with values in $\bS_{\delta}$, 
which is the space of $d\times d$ symmetric matrices  whose eigenvalues are in $[\delta,\delta^{-1}]$, and a Borel $\bR^{d}$-valued $b$ given on $\bR^{d}$.  

Fix some numbers 
$$
\rho_{a},\rho_{b}
\in(0,1],\quad p_{b}\in(1,\infty)
$$ and set
$$
\osc (a,B_{\rho}(x))=
 \dashint_{ B_{\rho}(x)}\dashint_{ B_{\rho}(x)}|a( y)-a( z)|\,dydz,
$$
$$
a^{\# }_{\rho_{a}}=\sup_{ x \in\bR^{d }}\sup_{\rho\leq \rho_{a}}
\osc (a,B_{\rho}(x)) ,
\quad
\hat b_{\rho_{b}}=\sup_{\rho\leq \rho_{b}}\rho \sup_{B\in\bB_{\rho}}\Big(\dashint_{B}|b|^{p_{b}}\,dx\Big)^{1/p_{b}},
$$
where $B_{\rho}(x)=\{y:|y-x|<\rho\}$,
$\bB_{\rho}$ the collection of $B_{\rho}(x)$,
$$
\dashint_{\Gamma}f\,dx=\frac{1}{|\Gamma|}
\int_{\Gamma}f\,dx
$$
and $|\Gamma|$ is the volume of $\Gamma$.

Recall that if $a^{\# }_{\infty}<\infty$, then one says that $a$
is in BMO (bounded mean oscillation) and, if $a^{\# }_{\rho_{a}}
\to 0$ as $\rho_{a}\to0$, then
$a\in VMO$ (vanishing mean oscillation). 

\begin{remark}
                          \label{remark 9.11.1}
1. Examples of $a\in VMO$
are given by (discontinuous functions) $\zeta(|x|)\ln|\ln|x||$,
$\zeta(|x|)\sin\ln|\ln|x||$,
where $\zeta$ is smooth and vanishes
for $|x|>1$. Examples of $b$ with finite $\hat b_{\rho_{b}}$ include $|b|=1/|x|$ with any $p_{b}\in(1,d)$.

2. If $b\in L_{d}$, then $\hat b_{\rho_{b}}$ for $p_{b}=d$
can be made as small as we like on account of 
choosing $\rho_{b}$ sufficiently small. This    is due to the fact that for $p_{b}=d$ we have
$$
\hat b_{\rho_{b}}=N(d)\sup_{B\in\bB_{\rho_{b}}}
\int_{B}|b|^{d}\,dx.
$$
\end{remark}

Our goal is
to show that under appropriate conditions
(on smallness of $a^{\# }_{\rho_{a}}$ and 
$\hat b_{\rho_{b}}$) on some 
complete probability spaces $(\Omega,\cF,P)$, provided with a filtration
of complete $\sigma$-fields $\cF_{t}
\subset \cF$, $t\geq0$, and a $d$-dimensional
Wiener (relative to $\{\cF_{t}\}$)
process $w_{t}$, $t\geq0$,
the equation
\begin{equation}
                                        \label{11.29.20}
x _{s}=x  +\int_{0}^{s}\sqrt{a(x_{r})}\,dw_{r}
+\int_{0}^{s}b ( x_{r}) \,dr
\end{equation}
has a solution and this solution
has unique finite-dimensional distributions (meaning that it is weakly unique).

We are interested in the so-called weak solutions, that is
solutions that are not necessarily $\cF^{w}_{s}$-measurable,
where $\cF^{w}_{s}$ is the completion of $\sigma(w_{u}:u\leq s)$.

After the classical work by K. It\^o showing that there exists
a unique (strong) solution of \eqref{11.29.20} if $\sigma$ and $b$
are Lipschitz continuous in $x$ (may also depend on $t$ and  $\omega$ and the nondegeneracy of $\sigma$ is not required),
  much effort was  applied to relax these Lipschitz conditions.   The first author who achieved a considerable progress was A.V. Skorokhod
\cite{Sk_61} who proved the solvability assuming  
only the continuity of $\sigma$ and $b$ 
in $x$ (which may depend
on $t$ and again without nondegeneracy).
Then by using the Skorokhod method and Aleksandrov
estimates the author proved in \cite{Kr_69}
and \cite{Kr_77} the solvability 
for the case of {\em just measurable $\sigma$ and  
bounded $b$}  under the nondegeneracy assumption.
Stroock and Varadhan \cite{SV_79} among other things
not only proved
the solvability for the coefficients uniformly 
continuous in $x$ (depending also on $t$), but also proved the uniqueness
of their distributions
(so called weak uniqueness).

M. R\"ockner and Guohuan  Zhao \cite{RZ_20} considered time-dependent coefficients with $a=(\delta^{ij})$ and, if restricted to our case, proved weak uniqueness (in a subsequent paper even strong uniqueness)
in case $b\in L_{d}$  (cf. Remark
\ref{remark 9.11.1}).

If $a=(\delta^{ij})$ the conditions for weak uniqueness can be significantly
reduced
further, see, for instance, \cite{Ki_24},
which also contains an excellent review
of previous results in this case. 

By the way,
G. Zhao (\cite{Zh_20_1}) gave an example showing
that, if in the definition of $\hat b_{\rho_{b}}$   we replace $\rho$ 
with $\rho^{\alpha}$, $\alpha>1$,
and require the new $\hat b_{\rho_{b}}<\infty$,   weak uniqueness
may fail even in the   case
of unit diffusion. 

 The  issue of weak uniqueness is also very well known to be related
to the operator 
  $$
\cL u( x)=(1/2) a^{ij}( x)D_{ij} u  ( x)+
b^{i}( x)D_{i} u ( x),
$$
acting on functions given on $\bR^{d}$
and to the solvability of the corresponding equations
in various function spaces.
Here and elsewhere the summation convention
over repeated indices is enforced. We also
use the notation $D_{i}u=\partial u/\partial x^{i}$, $Du=(D_{1}u,...,D_{d}u)$, $D_{ij}
=D_{i}D_{j}$,  $D^{2}u=(D_{ij}u)$.
\begin{remark}
                      \label{remark 8,9.4}
The following discussion is aimed at
 the readers accustomed to treating
elliptic equations in Sobolev spaces.

If $p_{b}>d$ the Sobolev space theory 
of solvability of elliptic and parabolic equations related to $\cL$ is
well known and is developed by using the 
perturbation method.

If $p_{b}\leq d$ and $b$ is a general element of $L_{p_{b}}$, then the operator
$ \cL$ seems to be  not even
continuous as an operator from $W^{ 2}_{p}$
into $L_{p} $. Indeed, if $p>p_{b}$,
the term $b^{i}D_{i}u$ will be generally 
not in $L_{p} $ because it may happen that $b\not\in L_{p}$. In case
$p=p_{b}$ for $b^{i}D_{i}u$ to be in $L_{p}$
it seems we have to have $Du$ bounded for any
$u\in W^{2}_{p}$ which is only the case
when $p>d$. In case $p<p_{b}$ for $b^{i}D_{i}u$ to be in $L_{p}$ H\"older's
inequality seems to require $Du\in L_{r}$
with $r=pp_{b}/(p_{b}-p)$ which by embedding
theorem is (only?) possible if $p_{b}\geq d$. 

Thus the only possible values for $p_{b}\leq d$ and $p$ for which the {\em Sobolev space\/} theory
seems to work is $p_{b}=d$ and $p<d$.
This case was treated in
\cite{Kr_21b}.  

The elliptic equations with $b$ having stronger singularities {\em expressed
in terms of $L_{p}$-spaces\/}
seem to not allow treatment in   Sobolev spaces. However, in case $b\in E_{p_{b},1}$ and $d\geq p_{b}
>  \sfd $,  where $\sfd=\sfd(d,\delta)\in (d/2,2)$, in
 \cite{Kr_23} there is a theory of
solvability of elliptic equations
in {\em Sobolev spaces\/} for $p<p_{b}$. Sending $p_{b}$ further down 
below $d/2$ seems to require the solvability theory
in Morrey spaces.

This is why we need Morrey spaces and Theorem \ref{theorem 10.27.1}.
Theorem \ref{theorem 6,6,1} will be used in a subsequent paper
to prove the strong Feller property and other properties
of the semigroup corresponding to the Markov process
constructed from solution of \eqref{11.29.20}.

\end{remark}

\mysection
{Parabolic Morrey spaces
and main results}

              \label{section 3.27.1}

Take $p,q\in(1,\infty)$. Define $L_{p}=L_{p}(\bR^{d})$, $L_{p,q}=L_{q}(\bR,L_{p})$, and let $L_{p,q,\loc}$ be the collection of functions
$f$ such that for any bounded set $A\subset
\bR^{d+1}=\{(t,x):t\in\bR,x\in\bR^{d}\}$,
we have
$fI_{A}\in L_{p,q}$. Let $C_{R}=[0,R^{2})\times B_{R}$,
$C_{R}(t,x)=(t,x)+C_{R}$ and let $\bC_{R}$
be the collection of $C_{R}(t,x)$, $(t,x)\in\bR^{d+1}$.

For    $\beta\geq 0$, introduce the
Morrey space $E_{p,q,\beta} $
as the set of $g\in  L_{p,q,\loc}$ such that  
\begin{equation}
                             \label{8.11.02}
\|g\|_{E_{p,q,\beta} }:=
\sup_{\rho\leq 1,C\in\bC_{\rho}}\rho^{\beta}
\dashnorm g  \|_{ L_{p,q}(C)} <\infty .
\end{equation}  
Define $\partial_{t}u=\partial u/\partial t$,
$$
E^{1,2}_{p,q,\beta} =\{u:u,Du,D^{2}u,
\partial_{t}u\in E_{p,q,\beta} \},
$$
and provide $E^{1,2}_{p,q,\beta} $ with an obvious norm.
The subsets of these spaces
consisting of functions independent of
$t$ is denoted by $E_{p,\beta}$
and $E^{2}_{p,\beta}$, respectively.
If $\cO$ is a subdomain of $\bR^{d+1}$
we denote $E_{p,q,\beta}(\cO)=\{u:u_{I_{\cO}}\in
E_{p,q,\beta}\}$,
$$ 
E^{1,2}_{p,q,\beta}(\cO) =\{u: u, Du, D^{2}u, 
\partial_{t}u\in E_{p,q,\beta}(\cO)\},
$$
and provide these spaces with  natural norms.

\begin{definition}
            \label{definition 7,16.1}
We call a     solution $x_{\cdot}$ of
 \eqref{11.29.20} with $x=0$ 
  {\em $E_{p,q,\beta }$-admissible\/} if
 for any $R\in(0,\infty)$
and Borel nonnegative $f$ on $\bR^{d+1}$
\begin{equation}
                      \label{12.11.06}
E\int_{0}^{R^{2}\wedge\tau_{R}}f(s,x_{s})\,ds\leq N\|f\|_{E_{p,q,\beta }},
\end{equation}
where $\tau_{R}$ is the first exit time of $ x_{s} $ from $B_{R}=B_{R}(0)$ and $N$ is independent of $f$. We call $x_{\cdot}$
  {\em $E_{p,\beta }$-admissible\/} if
 for any $R,t\in(0,\infty)$
and Borel nonnegative $f$ on $\bR^{d}$
\begin{equation}
                      \label{7,6.1}
E\int_{0}^{t\wedge\tau_{R}} f( x_{s})\,ds\leq N\|f\|_{E_{p,\beta }},
\end{equation}
where $N$ is independent of $f$.
Similarly we define admissible solutions
for any starting point $x$ (by taking it as the new origin).

\end{definition}

\begin{remark}
                    \label{remark 9.11.6}
Obviously if $x_{\cdot}$ is
$E_{p,q,\beta }$-admissible, then
it is also $E_{p,\beta }$-admissible.
\end{remark}

\begin{remark} 
                      \label{remark 9,2.1}
Owing to Theorem 1.5 of \cite{Kr_21}
any solution of \eqref{11.29.20}
with any starting point is $E_{\sfd,\beta}$-admissible for some $\sfd\in(d/2,d)$ and any
$\beta>0$ if $b\in L_{d}$
(note $d$ and not $\sfd$) since then the left-hand side of \eqref{7,6.1} is dominated by
$N\|f\|_{L_{\sfd}}\leq N\|f\|_{E_{\sfd,\beta}}$.

By Theorem 1.1 of \cite{Kr_21} the solutions
of \eqref{11.29.20} do exist if $b\in L_{d}$.

\end{remark}

To state our main results
fix  $p ,q ,\beta    $ such that
\begin{equation}
                        \label{3.21.01}
p ,q  \in(1,\infty),\quad\beta \in(1,\infty),\quad \beta\ne 2,  
\quad \frac{d}{p }+\frac{2}{q }\geq\beta .
\end{equation}

\begin{theorem}[weak uniqueness]  
              \label{theorem 7,6.1}
For $\beta<2$ and $p_{b}\geq p\beta$ there exist $\check a=\check a
(d,\delta,p,q)>0$ and $\check b=\check b_{0}
(d,\delta,p,q,\rho_{a},\beta)>0$ 
(defined in Section \ref{section 3.27.3}) such that, if $a^{\# }_{\rho_{a}}\leq \check a $ and $\hat b_{\rho_{b}}\leq \check b$, then
any two $E_{p,q,\beta}$- or $E_{p,\beta}$-admissible solutions  of
\eqref{11.29.20} with $x=0$ (if they exist) have the same finite-dimensional distributions.

\end{theorem}

\begin{remark}
On may be dissatisfied with the above ``conditional''
uniqueness of solutions. But let us recall that
the Lebesgue
uniqueness theorem for the one-dimensional
problem $u'=f$, $u(0)=0$, is also conditional:
you need $u$ to be absolutely continuous
and $f\in L_{1}$ and this {\em is not\/} ``only if'': say $u$ can be 
$x \sin(1/ x)$.
\end{remark}

The proof of Theorem \ref{theorem 7,6.1}
is given in Section \ref{section 3.27.3}.

\begin{remark}
                        \label{remark 9.11.2}
If $b\in L_{d}$, then solutions of \eqref{11.29.20} exist and they
are $E_{\sfd,\beta}$-admissible for any $\beta\geq0$. The condition regarding $\hat b_{\rho_{b}}\leq \check b$ is satisfied as well. But to have weak uniqueness we still need  $a^{\# }_{\rho_{a}}\leq \check a(d,\delta,\sfd,q)$, where $q\in(1,\infty)$ is arbitrary
(note that for $d\,(=p_{b})=\sfd\beta$
and any $q$ we have 
$\beta\leq d/\sfd+2/q$).
\end{remark}

\begin{remark}
                        \label{remark 9.11.4}
The smallness assumption of $a^{\# }_{\rho_{a}}$ for weak uniqueness is unavoidable as
shown in \cite{Na_97} and \cite{Sa_99}.

It turns out that without smallness assumption on $\hat b_{\rho_{b}}$ the equation even with unit diffusion may not
have {\em any\/} solution. For instance,
 let $a^{ij}=\delta^{ij}$,
$b(x)=-(d/2)I_{0<|x|\leq 1}x/|x|^{2}$
and assume that $x_{t}$ is a solution
of \eqref{11.29.20} with $x=0$.
Then
by It\^o's formula before $x_{t}$
exits from $B_{1}$ we have
$$
|x_{t}|^{2}=2\int_{0}^{t}x_{s}\,dw_{s}
+\int_{0}^{t}I_{x_{s}=0}\,ds,
$$ 
where the last integral is the time spent by $x_{s}$ at the origin before $t$. This time is zero, because
by It\^o's formula applied to
$|x^{1}_{t}|$ one easily sees that
$x^{1}_{t}$ has finite local time at 0, which implies that the real time
it spends at zero is zero. Then the above formula implies  that $|x_{t}|^{2}$ is a nonnegative local martingale
and therefore it equals zero. However,
$x_{t}\equiv 0$ does not satisfy~\eqref{11.29.20}.  
\end{remark}

Assumption that $\hat b_{\rho_{b}}\leq
\check b$
becomes stronger if in the definition of 
$\hat b$ we replace $p_{b}$ with any
$p\geq p_{b}$. It is instructive to know
that if such modification holds
with $p\in(\sfd,d)$ it might happen
that $b\not\in L_{p+\varepsilon,\loc}$,
no matter how small $\varepsilon>0$ is.
\begin{example}[see also Example
\ref{example 7,29.1}]
                        \label{example 9.25.1}

 Take $p \in[d-1,d)$ and take $r_{n}>0$, $n=1,2,...$, such that
the  sum of $\rho_{n}:=r_{n}^{d-p}$ is $1/2$. Let $e_{1}$ be the first
basis vector  and set $b(x)=|x|^{-1}
I_{|x|<1}$, 
$b_{n}(x)=r_{n}^{-1}b\big(r_{n}^{-1}
(x-c_{n}e_{1})\big)$, $n=1,2,...$, $\tilde b=\sum _{n} b_{n}$,  where
$x_{0}=1$,
$$
x_{n}=1-  2\sum_{1}^{n}r_{i}^{d-p},\quad 
c_{n}=(1/2)(x_{n}+x_{n-1}).
$$

Since $r_{n}\leq 1$ and $d-p\leq 1$,  the supports of $b_{n}$'s are disjoint and
for $q>0$
$$
\int_{B_{1}}\tilde b^{q}\,dx=\sum _{1}^{\infty}\int_{\bR^{d}}b_{n}^{q}\,dx=N(d,p)\sum_{1}^{\infty}r_{n}^{d-q}.
$$
According to this we take the $r_{n}$'s so that
the last sum diverges for any $q>p$.
Then observe that for any $n\geq 1$ and any ball $B$
of radius $\rho$
$$
 \int_{B }  b_{n} ^{p}dx \leq N(d) \rho^{d-p} .
$$
Also, if the intersection of $B$ with $\bigcup B_{r_{n}}(c_{n})$
is nonempty, the intersection
 consists of some $B_{r_{i}}(c_{i})$, $i=i_{0},...,i_{1}$, and $B\cap B_{r_{i_{0}-1}}(c_{i_{0}-1})$ if $i_{0}\geq 2  $ and 
$B\cap B_{r_{i_{1}+1}}(c_{i_{1}+1})$.
In this situation $ c_{i_{0}}-c_{i_{1}}
\leq 2\rho$, $c_{i_{0}}\leq x_{i_{0}},
c_{i_{1}}\geq x_{i_{i-1}}$ and
$$
 \sum_{i=i_{0} +1}^{i_{1}-1 }
r^{d-p}_{i} \leq    \rho.
$$
Therefore,
$$ \int_{B }\tilde  b  ^{p}\,dx\leq N(d)\sum_{i=i_{0}+1}^{i_{1}-1}r_{i}^{d-p}
+\int_{B }  [I_{i_{0}\geq2} b_{i_{0}-1}  ^{p} 
+b^{p}_{i_{0}}+ b_{i_{1} }  ^{p}+ b_{i_{1}+1}  ^{p}]dx \leq N(d)(\rho
+\rho^{d-p}),
$$
where the last term is less than $N(d)\rho^{d-q}$
for $\rho\leq 1$. This domination
also holds for $\rho>1$,
since $b\in L_{p}$ and $d\geq p$, and this yields  $\hat b_{\rho_{b}}\leq
\check b$ for $c\tilde b$ in place of $b$
with $c$ small enough.
\end{example}

As usual the proof  of
Theorem \ref{theorem 7,6.1} is based on a version of It\^o's formula.

\begin{theorem}[It\^o's formula]
                 \label{theorem 3.16.1}
Let \eqref{3.21.01} be satisfied, let $1<\beta<2$, $p_{b}\geq p\beta$, and $\hat b
_{\rho_{b}}<\infty$.

(i) If
  $x_{\cdot}$   is an $E_{p,q,\beta }$-admissible solution  of
 \eqref{11.29.20} with $x=0$ then
 for any $R\in(0,\infty)$ 
  and   $u\in E^{1,2}_{p,q,\beta  } $, 
 with probability one for all $t\in (0,R^{2}]$,
\begin{equation}
                                \label{3.16.2}u(t\wedge\tau_{R},x_{t\wedge\tau_{R}})
=u( 0)+\int_{0}^{t\wedge\tau_{R}}D_{i}u (s,x_{s})
\sigma^{ik}(x_{s}) \,dw^{k}_{s}
+\int_{0}^{t\wedge\tau_{R}}[
\partial_{t}u+\cL u] (s,x_{s})\,ds,
\end{equation}
where $\sigma=\sqrt a$, and the stochastic integral above is a 
square-integrable
martingale.

(ii) If
  $x_{\cdot}$   is an $E_{p, \beta }$-admissible solution  of
 \eqref{11.29.20} with $x=0$ then
  for any $R\in(0,\infty)$,
  and   $u\in E^{1,2}_{p,\beta  } $  
 with probability one for all $t\geq0$,
$$
u( x_{t\wedge\tau_{R}})
=u( 0)+\int_{0}^{t\wedge\tau_{R}}D_{i}u ( x_{s})
\sigma^{ik}(x_{s}) \,dw^{k}_{s}
+\int_{0}^{t\wedge\tau_{R}} 
 \cL u  ( x_{s})\,ds,
$$
 and the stochastic integral above is a 
square-integrable
martingale

\end{theorem}

This theorem is proved in Section 
\ref{section 9.12.1}.

\begin{remark}
                 \label{remark 7,6.1}
It is, certainly, desirable to express
conditions \eqref{12.11.06} and \eqref{7,6.1} in terms of usual $L_{p,q}$
and $L_{p}$ spaces instead of Morrey
spaces. Therefore, it is useful
to have in mind that, on the right-hand
side of  \eqref{12.11.06} one can replace
$f$ with $fI_{C_{R}}$, where $C_{R}=
[0,R^{2})\times B_{R}$, which shows that,
if \eqref{12.11.06} holds with
  $L_{p,q}(C_{R})$ in place
of $L_{p,q,\beta}$, then  it holds
as is. Similar observation 
is valid for \eqref{7,6.1}.

However, even if $x_{t}$ is a Wiener
process, so modified \eqref{12.11.06}
holds only when $d/p+2/q<2$. In case
of modified \eqref{7,6.1} we need
$d/p<2$. Therefore, the modified
conditions \eqref{12.11.06} and \eqref{7,6.1} become unrealistic if, say
$p$ is small. At the same time 
one can show that  under smallness conditions
on $a^{\# }_{\rho_{a}}$ and $\hat b_{\rho_{b}}$ (including the case
of the Wiener process without drift)
{\em there are solutions\/} satisfying
  \eqref{12.11.06} and \eqref{7,6.1}   as long
as $1<\beta\leq d/p+2/q$ in the case
of \eqref{12.11.06} and $1<\beta\leq d/p$ in the case of \eqref{7,6.1}.

\end{remark}

\begin{remark}
                 \label{remark 9.12.1}
We know that  
under our assumptions $E_{p,\beta}$-admissible solutions
of \eqref{11.29.20}  exist. However, the page restriction does not allow us to include
the proof of it. We can only say that
it is quite different from  the proof
given in \cite{1} where $p_{b}>d/2$
and the solvability of   equations in Sobolev spaces is used. In case $p_{b}\leq d/2$
our proof instead is based on Theorems \ref{theorem 6,6,1} and \ref{theorem 10.27.1}.

\end{remark}

\mysection{Properties of 
$E_{p,q,\beta}$}

                \label{section 9.12.1}

We are going to use some results
 proved in \cite{Kr_22} and \cite{2}. 
Observe that in \cite{Kr_22} the $E_{p,q,\beta}$-norm
is defined differently form \eqref{8.11.02}
without the restriction on $\rho$, just $\rho<\infty$.
In connection with this notice that if $\beta\leq
d/p+2/q$ and the support of $g$ is in some $C\in\bC_{1}$,
then the $E_{p,q,\beta}$-norm of $g$ is the same:
taken from \eqref{8.11.02} or from \cite{Kr_22}.
One more useful observation is that
\begin{equation}
                               \label{3.30.3}
\|g\|_{E_{p,q,\beta} }=\sup_{C\in \bC_{1}}
\|gI_{C}\|_{E_{p,q,\beta} }.
\end{equation}

A very unusual property of $E_{p,\beta}$
is that smooth functions are not dense 
there. Even more, the distance in
$E_{d/2,1}$ between $u=|x|^{-1}$ and any smooth function is the same as its distance to zero.

However, here is a useful approximation result. 
For functions $f(t,x)$ and $\varepsilon>0$ we define
\begin{equation}
                        \label{7,5.1}
f^{(\varepsilon)}=f*\zeta_{\varepsilon},
\end{equation}
 where
  $\zeta_{\varepsilon}=\varepsilon^{-d-2}\zeta
(t/\varepsilon^{2},x/\varepsilon)$ with a nonnegative
$\zeta\in C^{\infty}_{0}(\bR^{d+1})$ which has unit integral and support
in $C_{1}(-1,0)$. Observe that
owing to Minkowski's inequality  
$$
\|f^{(\varepsilon)}\|_{E_{p,q,\beta}}
\leq \|f \|_{E_{p,q,\beta}}
$$
for any $f\in E_{p,q,\beta}$.
\begin{lemma}[Lemma 2.3 of \cite{2}]
                               \label{lemma 3.14.3}
Let $0\leq\beta'<\beta$ and 
  $f\in E_{p,q,\beta'}$, then for any $C\in \bC$
\begin{equation} 
                            \label{3.14.10}
\lim_{\varepsilon\downarrow0}
\|f^{(\varepsilon)}-f \|_{E_{p,q,\beta}(C) }=0.
\end{equation}
\end{lemma}

Another result of approximation type
is the following.    
\begin{lemma}
                \label{lemma 7,4.3}
Let   $g(t,x)\geq 0$ be a Borel function such that for any smooth bounded $f(t,x)$ we have
\begin{equation}
                       \label{7,4.5}
\int_{\bR^{d+1}}g|f|\,dxdt\leq
\|f\|_{E_{p,q,\beta}}.
\end{equation}
Then, for any $f\in E_{p,q,\beta}$,
\eqref{7,4.5} holds and, moreover,
$$
\lim_{\varepsilon\downarrow 0}
\int_{\bR^{d+1}}
g|f-f^{(\varepsilon)}|\,dxdt=0.
$$
\end{lemma}

 \begin{remark} 
                        \label{remark 7,2.1}
Holder's inequality easily shows that
$$
\|f\|_{E_{p,\beta}}\leq N
\sup_{B\in \bB_{1}}\|f\|_{L_{r}(B)},
$$
where $N$ is independent of $f$ and 
$r=p\vee(d/\beta)$. Also,
if $\beta\leq d/p+2/q$, then
$\|g\|_{E_{p,q,\beta}}\leq \|g\|_{L_{\bar p,\bar q}}$, where $(\bar p,\bar q) =\alpha (p,q)$ and $\beta\alpha=
d/p+2/q$.
\end{remark}

\begin{remark}
                        \label{remark 3.27.10}
If $u$ vanishes outside $C\in\bC_{1}$, then, if $d/p +2/q \geq\beta $, one easily sees that
for any $r\geq 1$ and $C'\in \bC_{r}$
$$
r^{\beta }\dashnorm u\|_{L_{p,q}(C')}\leq N(d)
\|u\|_{L_{p,q}(C)}.
$$
\end{remark}

Remark 5.8 of \cite{Kr_22} imply
the following easy consequence
of H\"older's inequality.

\begin{lemma}
                \label{lemma 6,11.1}
If $p,q\in(1,\infty]$,
$d/p+2/q\geq \beta >1$, and

$(p_{0},q_{0})=\beta(p,q)=(s,r)(\beta-1) $, then for any
$f,g$
$$
\|fg\|_{E_{p,q,\beta}}
\leq \|f\|_{E_{p_{0},q_{0},1}}\|g\|_{E_{s,r,\beta-1}}.
$$

\end{lemma}

We have the following as a consequence of the interpolation Lemma  5.10 of 
\cite{Kr_22}. 

\begin{lemma}
                        \label{lemma 3.27.4}
If $0<\beta\leq d/p+2/q$, then 
for any $u\in E^{1,2}_{p,q,\beta}$ and $\varepsilon\in(0,1]$,
\begin{equation}
                                    \label{3.20.6}
\| Du\|_{E_{p,q,\beta}}\leq \varepsilon 
\|\partial_{t}u,D^{2}u\|_{E_{p,q,\beta}}
+N(d,p,q,\beta)\varepsilon^{-1}\| u\|_{E_{p,q,\beta}}.
\end{equation} 
\end{lemma}

If $u$ is independent of $t$, \eqref{3.20.6}
provides and interpolation inequality in the ``elliptic"
case.

The following embedding lemma is 
part of Lemma 2.5 of \cite{2}.
As usual, when we say that a function
$u$
uniquely defined only up to almost everywhere is, say continuous, we mean
that there is a continuous function
equal to $u$ almost everywhere.
 
\begin{lemma}
           \label{lemma 3.27.40} 
  Let   $ \beta <2 $. Then any $u\in E^{1,2}_{p,q,\beta}$
is bounded and continuous (meaning that it has a modification which is 
bounded and continuous)
 and for any $\varepsilon
\in(0,1]$  
\begin{equation}
                                    \label{3.20.06}
|u|
\leq \varepsilon^{2-\beta} \|\partial_{t}u,D^{2}u\|_{E_{p,q,\beta}}+N(d, \beta)\varepsilon^{-\beta} \| u\|_{E _{p,q,\beta}}
\end{equation} 
with $\beta=2-d/p-2/q$ in the powers of $\varepsilon$.
 
\end{lemma}

Another useful property of $E^{1,2}_{p,q,\beta}$
is an embedding theorem, which
follows from Corollary 5.7 of \cite{Kr_22}.

\begin{lemma}
                \label{lemma 3.27.5}
Let $ 1< \beta\leq d/p+2/q$. Then for any $u\in E^{1,2}_{p,q,\beta}$
we have
\begin{equation}
                                    \label{3.27.5}
\|Du\|_{E_{r,s,\beta-1}}\leq N(d,p,q,\beta)  \| u\|_{E^{1,2}_{p,q,\beta}},
\end{equation} 
where $(r,s)(\beta-1)=(p,q)\beta$.
\end{lemma}
 
Lemmas \ref{lemma 6,11.1} and \ref{lemma 3.27.5} imply the following 
consequence of Lemma
5.1 of \cite{2}.   

\begin{theorem}
              \label{theorem 8,23.3}
Let $p,q\in(1,\infty)$,
$$
\frac{d}{p}+\frac{2}{q}\geq \beta>1,
$$ 
$\beta p\leq p_{b}$, and $\hat b_{\rho_{b}}<\infty$. Then the operator
$\cL$ is a bounded operator
from $E^{1,2}_{p,q,\beta}$ to $E_{p,q,\beta}$. 

\end{theorem}

{\bf Proof of Theorem \ref{theorem 3.16.1}}. 
Proof. (i) Observe that $d/p
+2/q\geq\beta>1
 $ and by Lemma \ref{lemma 3.27.5} we have
$|Du|^{2}\in E_{r/2,s/2,2(\beta -1) }$, where $r=
p\beta /(\beta -1)$,  $s=
q\beta  /(\beta -1)$. Note that  
$$
  2>\beta  >1,\quad \beta  /(\beta   -1) >2,
\quad 2(\beta  -1)  < \beta ,
\quad
r /2\geq p, \quad s /2\geq q.
$$
It follows that $E_{r/2,s/2,2(\beta  -1) }
\subset E_{p,q,\beta }$ and this along 
with \eqref{12.11.06}  implies  the last statement in (i).

The inequality $ \beta  -1   < \beta/2$
and Lemma \ref{lemma 3.14.3} also imply that
$$
E\Big|\int_{0}^{t\wedge\tau_{R }}D_{i}\big(u^{(\varepsilon)}-u\big) 
\sigma^{ik} (s,x_{s})\,dw^{k}_{s}\Big|^{2}
$$
$$
\leq N\int_{0}^{t\wedge\tau_{R }}|D  u^{(\varepsilon)}-Du|^{2}(s,x_{s})\,ds
\leq N\|D  u^{(\varepsilon)}-Du\|^{2}_{E_{2p,2q, \beta/2 }(C_{R\vee\sqrt t})}\to0
$$
as $\varepsilon\downarrow 0$.
This shows that after we apply It\^o's
formula to $u^{(\varepsilon)}$, we will be able to pass to the limit in the stochastic integral.

In what concerns the usual integral
observe that   estimate \eqref{12.11.06} implies that, for any $R$, there is a Borel
function $g(t,x)\geq0$ such that
$$
\int_{C_{R}}gf\,dxdt=E\int_{0}^{R^{2}\wedge\tau_{R}}f(s,x_{s})\,ds\leq N\|f\|_{E_{p,q ,\beta  }} 
$$
for any $f\geq0$.
Then Lemma \ref{lemma 7,4.3} implies that
$$
\int_{C_{R}}g|a^{ij}D_{ij}(u^{(\varepsilon)}-u)| \,dxdt\leq
N\int_{C_{R}}g| (D u)^{(\varepsilon)}-Du| \,dxdt\to 0
$$
as $\varepsilon\downarrow 0$.
This shows that we can pass to the limit in the usual integral containing
$a^{ij}D_{ij}u^{(\varepsilon)}$.
Furthermore, for any $R>0$ there is a Borel $h\geq0$ such that
$$
\int_{C_{R}}h f\,dxdt=E\int_{0}^{R^{2}\wedge\tau_{R}}|b|f(s,x_{s})\,ds\leq N\|bf\|_{E_{p,q,\beta }} 
\leq N\|f\|_{E_{r,s,\beta-1}}
$$
for any $f\geq0$, where
the last inequality follows from
Lemma  \ref{lemma 6,11.1}.
Then Lemma \ref{lemma 7,4.3} implies that
$$
E\int_{0}^{R^{2}\wedge\tau_{R}}
|b|  \,|D u^{(\varepsilon )}-Du |(s,x_{s})\,ds
=\int_{C_{R}}h |D u^{(\varepsilon )}-Du |\,dxdt\to 0
$$
as $\varepsilon\downarrow 0$.
Finally, since $u$ is bounded and continuous (Lemma \ref{lemma 3.27.40}), we have
the convergence  of the terms 
$u^{(\varepsilon)}(t\wedge\tau_{R},x_{t\wedge\tau_{R}})
,u^{(\varepsilon)}( 0)$ to
$u (t\wedge\tau_{R},x_{t\wedge\tau_{R}})
,u ( 0)$.
This proves (i).
To prove (ii),
it suffices to take the above $u$
independent of $t$.  The theorem is proved.  \qed

The following theorem, as we have mentioned in Remark \ref{remark 9.12.1},
is needed to show the solvability
of \eqref{11.29.20} if $p_{b}\leq d/2$. We show that the $t$-traces
of functions in $E^{1,2}_{p,q,\beta}$
possess some regularity as $L_{p}$-functions. For $\gamma=0$ or $1$ set
$$
D^{\gamma}=D\quad\text{if}\quad
\gamma=1\quad\text{and}\quad D^{\gamma}=1\quad
\text{if}\quad \gamma=0.
$$
Below by
$D^{\gamma}u(0,\cdot)$ we mean the limit
in $L_{r}(B)$ for any ball $B$
of $D^{\gamma}u^{(\varepsilon)}(0,\cdot)$
as $\varepsilon\downarrow 0$.
The existence of this limit easily follows from Lemma \ref{lemma 6,17.1} and Corollary \ref{corollary 6,19.1} below.
\begin{theorem}
              \label{theorem 6,6,1}
Take  
$r\in[ p,\infty)$, $\mu>0$    and assume
that  
$$
2-\gamma< \beta \leq \frac{d}{p}+\frac{2}{q}
<2-\gamma+\frac{d}{r},
\quad \kappa:=\gamma+\frac{d}{p}
+\frac{2}{q}-\frac{d}{r}\leq \mu
<2 .
$$

 Then for any
$u\in E_{p,q,\beta}$ the trace $D^{\gamma}u(0,\cdot)$ is uniquely defined
and for any $\varepsilon>0$
\begin{equation}
                     \label{6,6.1}
\|D^{\gamma}u(0,\cdot)\|_{E_{r,\beta+\gamma-\mu}(\bR^{d})}\leq N\varepsilon \|\partial_{t}u,
D^{2}u\|
_{ E _{p,q,\beta}}
+N\varepsilon^{-\mu/(2-\mu)}
\|u\|
_{ E _{p,q,\beta}},
\end{equation} 
\begin{equation}
                     \label{8,27.1}
\|D^{\gamma}u(0,\cdot)\|_{E_{r,\beta+\gamma-2}(\bR^{d})}\leq N  \|u\|
_{ E^{1,2} _{p,q,\beta}},
\end{equation} 
where  the constants $N$ depend only
on $d,p,q,\beta,\mu$.

\end{theorem}

To prove Theorem \ref{theorem 6,6,1},
first, we need the following corollary of
Theorem 10.2 of \cite{BIN_75} which
we give with a different proof for completeness.

\begin{lemma}
                 \label{lemma 6,17.1}
Let  
$r\geq p$  and assume
that $\kappa<2$.
 Then for any
$u\in W^{1,2}_{p,q}=\{u:\partial_{t}u,D^{2}u,Du,u\in L_{p,q}\}$ and $\varepsilon>0$ we 
have
\begin{equation}
                            \label{6,17.2}
\|D^{\gamma}u(0,\cdot)\|_{L_{r}}\leq
N\varepsilon \|\partial_{t}u,D^{2}u
\|_{L_{p,q} }+N\varepsilon^{-\kappa/(2-\kappa)}
\|u
\|_{L_{p,q} }.
\end{equation}
\end{lemma}

Proof. The case of arbitrary $\varepsilon>0$
is reduced to that of $\varepsilon=1$
by using   parabolic scalings.
To treat $\varepsilon=1$ take $\zeta\in
C^{\infty}_{0}(\bR)$ such that $\zeta(t)=1$
for $|t|\in[0,1]$, $\zeta(t)=0$ for $|t|\geq 2$, and define $-f=
\partial_{t}(\zeta u)+\Delta(\zeta u)$. Introduce
$$
P_{\gamma}(t,x)=t^{-(d+\gamma)/2}e^{-|x|^{2}/(8t)}.
$$
We know (from It\^o's formula or
from PDEs) that $\zeta u=Rf$, where
$$
Rf(t,x)=N(d)\int_{0}^{\infty}s^{-d/2}
\int_{\bR^{d}}e^{-|x-y|^{2}/(4s)}
f(t+s,y)\,dyds.
$$
It follows that ($y^{0}:=1,y^{1}:=y$)
$$
D^{\gamma}u(0,x)= -N(d)\int_{0}^{\infty}t^{-d/2}
\int_{\bR^{d}}\Big(\frac{y}{2t}\Big)^{\gamma}e^{-| y|^{2}/(4t)}f(t,x-y)\,dydt.
$$
By observing that $|y/\sqrt t|
e^{-| y|^{2}/(4t)}\leq Ne^{-| y|^{2}/(8t)}$ we conclude that
$
|D^{\gamma}u(0,x)|\leq NF(x),
$
where
$$
 F(x) =
\int_{0}^{2} 
\int_{\bR^{d}}P_{\gamma}(t,y)|f(t,x-y)|\,dydt .
$$

By Minkowski's inequality
$$
\|F\|_{L_{r} }
\leq \int_{0}^{2}\Big(
\int_{\bR^{d}}\Big(\int_{\bR^{d}}
P_{\gamma}(t,x-y)|f(t,y)|\,dy\Big)^{r}
\,dx\Big)^{1/r}
\,dt,
$$ 
where inside the integral with respect
to $t$ we have the norm of a convolution, so that by Young's
inequality this expression is dominated by
$
\|f(t,\cdot)\|_{L_{p}}\|P_{\gamma}(t,\cdot)
\|_{L_{s}},
$
where $1/s=1+1/r-1/p$ ($ \leq 1$ since $r\geq p$). An easy computation
shows that $\|P_{\gamma}(t,\cdot)
\|_{L_{s}}=N(d)t^{\alpha}$ with
$\alpha=-\gamma/2+(d/2)(1/r-1/p)$, which  yields
$$
\|F\|_{L_{r} }
\leq N\int_{0}^{2}\|f(t,\cdot)\|_{L_{p} }t^{\alpha}\,dt.
$$
Now use H\"older's inequality along with the observation that
$\alpha q/(q-1)>-1$, due to 
the assumption that $\kappa<2$, that is,  $2-\gamma +d/r>
d/p+2/q$, to conclude that
$
\|F\|_{L_{r} }
\leq N \|f \|_{L_{p,q}  }.
$ 
The lemma is proved.\qed

\begin{corollary}
           \label{corollary 6,19.1}
For any $\rho\leq 1,\varepsilon>0$ and
$u\in W^{1,2}_{p,q}(C_{2\rho})$ we have
$$
\dashnorm D^{\gamma} u(0,\cdot)\|_{L_{r}(B_{\rho})}
\leq N\varepsilon  
\dashnorm \partial_{t}u,D^{2}u
\|_{L_{p,q}(C_{2\rho})}
$$
$$
+N
(\varepsilon \rho^{-2}+\varepsilon^{-\kappa/(2-\kappa)}\rho^{(2\kappa-2\gamma)/(2-\kappa)})
\dashnorm u
\|_{L_{p,q}(C_{2\rho})}.
$$
\end{corollary}

Indeed, the case of $\rho<1$ is reduced
to $\rho=1$ by means of parabolic dilation. In the latter case
take $\zeta\in C^{\infty}_{0}
(\bR^{d+1})$ such that $\zeta=1 $ on $C_{1}$ and $\zeta=0$ in $\bR^{d+1}_{0}
\cap C_{2}$. Then  use
\eqref{6,17.2} to see that
$$
\|D^{\gamma}u(0,\cdot)\|_{L_{r}(B_{1})}
\leq N\varepsilon  
\|\partial_{t}(\zeta u),D^{2}(\zeta u)
\|_{L_{p,q} }+N\varepsilon^{-\kappa/(2-\kappa)}
\|u
\|_{L_{p,q}(C_{2}) }
$$
$$
\leq N\varepsilon  
\|\partial_{t}u,D^{2}u
\|_{L_{p,q}(C_{2}) }+N\varepsilon \|u,Du\|_{L_{p,q}(C_{2}) }
+N\varepsilon^{-\kappa/(2-\kappa)}
\|u
\|_{L_{p,q}(C_{2}) }.
$$
After that it only remains to use
the interpolation inequality
$$
\|Du(t,\cdot)\|_{L_{p}(B_{2})}\leq
N \|D^{2}u(t,\cdot)\|_{L_{p}(B_{2})}+
N\| u(t,\cdot)\|_{L_{p}(B_{2})}.
$$

To prove  Theorem \ref{theorem 6,6,1} we also need its homogeneous 
version for homogeneous Morrey
space $\dot E_{p,q,\beta} $
introduced in the same way as $ E_{p,q,\beta} $ with the only difference that $\rho\leq 1$
in \eqref{8.11.02} is replaced with
$\rho<\infty$.
Accordingly define
$$
\dot E^{1,2}_{p,q,\beta} =\{u:u,Du,D^{2}u,
\partial_{t}u\in  \dot E_{p,q,\beta} \},
$$
and provide $\dot E^{1,2}_{p,q,\beta} $ with an obvious norm.

\begin{lemma}
               \label{lemma 6,18.1}
Let   
$r\geq p$  and let
$$
2-\gamma< \beta \leq \frac{d}{p}+\frac{2}{q}
<2-\gamma+\frac{d}{r}.
$$ 
Then for any
$u\in \dot E^{1,2}_{p,q,\beta}$ its trace $u(0,\cdot)$ is uniquely defined
and  
\begin{equation}
                     \label{6.6.1}
\|D^{\gamma}u(0,\cdot)\|_{\dot E_{r,\beta+\gamma-2}(\bR^{d})}\leq N \|\partial_{t}u,
D^{2}u\|
_{ \dot E _{p,q,\beta}},
\end{equation} 
where the constant  $N$ depends only
on $d,p,q,r$.

\end{lemma}

Proof. Take $\zeta\in C^{\infty}_{0} (\bR^{d+1})$, such that $\zeta(0)=1\geq \zeta\geq0$, define $\zeta_{n}(t,x )=
\zeta(t/n^{2},x/n )$ and observe that,
as $n\to\infty$,
$$
 \big|\|
\partial_{t}(\zeta_{n}u)\|_{\dot E_{p,q,\beta}} -\|\zeta_{n}
\partial_{t}u\|_{\dot E_{p,q,\beta} }\big|\leq n^{-2}\sup |\partial_{t}\zeta|
\|
u\|_{\dot E_{p,q,\beta}} \to 0.
$$
Also
$$
 \big|\|
D(\zeta_{n}u)\|_{\dot E_{p,q,\beta}} -\|\zeta_{n}
Du\|_{\dot E_{p,q,\beta} }\big|\leq n^{-1}\sup |D\zeta|
\|
u\|_{\dot E_{p,q,\beta}} \to 0,
$$
$$
 \big|\|
D^{2}(\zeta_{n}u)\|_{\dot E_{p,q,\beta}} -\|\zeta_{n}
D^{2}u\|_{\dot E_{p,q,\beta} }\big|\leq n^{-2}\sup |D^{2}\zeta|
\|
u\|_{\dot E_{p,q,\beta}}
$$
$$
+2n^{-1}\sup |D\zeta|\|
Du\|_{\dot E_{p,q,\beta} } \to 0.
$$

It follows that it suffices to concentrate on $u$ that vanish
for large $|t|+|x|^{2}$. In that case set $-f=\partial_{t}u+\Delta u$.
To further reduce our
problem observe that using translations show that it suffices
to prove that for any
  $0<\rho\leq 1$,
\begin{equation}
                     \label{6.19.1}
\rho^{\beta+2-\gamma}\dashnorm D^{\gamma}u(0,\cdot)\|_{L_{r}( B_{\rho})}\leq N\sup_{\rho_{1}\geq \rho}
\rho^{\beta }_{1}\dashnorm f\|
_{L_{p,q} (C_{\rho_{1}})   )}
=  N\sup_{\rho_{1}\in[\rho,\rho_{2}]}
\rho^{\beta }_{1}\dashnorm f\|
_{L_{p,q} (C_{\rho_{1}})   },
\end{equation}
where $\rho_{2}$ is such that
$u(t,x)=0$ for $|t|+|x|^{2}\geq \rho_{2}^{2}$
and the last equality is due to
$\beta \leq d/p+2/q$.

It is easy to pass to the limit
in \eqref{6.19.1} from smooth functions to arbitrary ones in
$W^{1,2}_{p,q}(C_{\rho_{2}})\supset
E^{1,2}_{p,q,\beta}(C_{\rho_{2}})$.
Therefore, we may assume that $u$ is smooth.
We thus reduced the general case to the task of proving the first estimate
in  \eqref{6.19.1} for smooth $u$ with compact support. One more reduction
is achieved by using parabolic scalings
which show that we only need to concentrate on $\rho=1$, that is
we only need to prove
\begin{equation}
                    \label{6,19,1}
 \dashnorm D^{\gamma}u(0,\cdot)\|_{L_{r}( B_{1})}\leq N\sup_{\rho \geq 1}
\rho^{\beta} \dashnorm f\|
_{L_{p,q} (C_{\rho })   )}
\end{equation}
for smooth $u$ with compact  support.

  Now   define
  $
g=|f|I_{C_{2}},h=|f|I_{C_{2}^{c}}$.
As it follows from the proof of Lemma~
\ref{lemma 6,17.1},
$$
|D^{\gamma}u(0,x)|\leq NG_{\gamma}(x)+NH_{\gamma}(x),
$$
where
$$
 (G_{\gamma},H_{\gamma})(x) =
\int_{0}^{\infty} 
\int_{\bR^{d}}P_{\gamma}(t,x-y)(g,h)(t,y)\,dydt.
$$

Estimate (2.2) of \cite{Kr_22}
(which allows $\alpha=0$ there)
implies that
that for $|x|\leq 1$
$$
H_{\gamma}(x)\leq N\sup_{\rho>1}\rho^{\beta}
\dashint_{(0,x)+C_{\rho}}h\,dydt
\leq N\sup_{\rho>1}\rho^{\beta}
\dashint_{ C_{2\rho}}h\,dydt
$$
$$
\leq N\sup_{\rho \geq1}
\rho^{\beta} \dashnorm f\|
_{L_{p,q}  (C_{\rho} )},
$$
where the last inequality is
due to H\"older's inequality.
Hence,
\begin{equation}
                     \label{6,8,3}
 \| H_{\gamma}\|_{L_{r}(  B_{1})}\leq N\sup_{\rho \geq1}
\rho^{\beta} \dashnorm f\|
_{L_{p,q}  (C_{\rho} )}. 
\end{equation}

By Minkowski's inequality
$$
\|G_{\gamma}\|_{L_{r}(B_{1})}
\leq \int_{0}^{\infty}\Big(
\int_{B_{1}}\Big(\int_{\bR^{d}}
P_{\gamma}(t,x-y)g(t,y)\,dy\Big)^{r}\,dx\Big)^{1/r}
\,dt,
$$ 
where inside the integral with respect
to $t$ we have the norm of a convolution, so that by Young's
inequality this expression is dominated by
$$
\|g(t,\cdot)\|_{L_{p}}\|P_{\gamma}(t,\cdot)
\|_{L_{s}},
$$
where $1/s=1+1/r-1/p$ ($ \leq 1$ since $r\geq p$). We know that $\|P(t,\cdot)
\|_{L_{s}}=N(d)t^{\alpha}$ with
$\alpha=-\gamma/2+(d/2)(1/r-1/p)$, which after taking into account that 
$\|g(t,\cdot)\|_{L_{p}}=0$ for $t\geq4$, yields
$$
\|G\|_{L_{r}(B_{1})}
\leq N\int_{0}^{4}\|f(t,\cdot)\|_{L_{p}(B_{2})}t^{\alpha}\,dt.
$$
Now use H\"older's inequality along with the observation that
$\alpha q/(q-1)>-1$, due to 
the assumption that   $2-\gamma+d/r>
d/p+2/q$, to conclude that
$$
\|G\|_{L_{r}(B_{1})}
\leq N \|f \|_{L_{p,q}(C_{2})}.
$$
This and \eqref{6,8,3} prove
\eqref{6,19,1} and the lemma. \qed

{\bf Proof of Theorem \ref{theorem 6,6,1}}.  To prove \eqref{6,6.1}, it suffices to show that
for any $\rho\in (0, 1]$,
$\varepsilon>0$ 
\begin{equation}
                     \label{6,20.1}
I_{\rho}:=\rho^{\beta+\gamma-\mu}\dashnorm Du(0,\cdot)\|_{L_{r}(B_{\rho})}\leq N\varepsilon \|\partial_{t}u,
D^{2}u\|
_{ E _{p,q,\beta}}
+N\varepsilon^{-\mu/(2-\mu)} 
\|u\|
_{ E _{p,q,\beta}}.
\end{equation} 

By Corollary
\ref{corollary 6,19.1} with
$\epsilon=\varepsilon \rho^{ \gamma-\mu }$
in place of $\varepsilon$  we get
$$
 I_{\rho}\leq N\epsilon \|\partial_{t}u,
D^{2}u\|
_{ E _{p,q,\beta}}+N\big(\epsilon  \rho^{-2} +\epsilon^{-\kappa/(2-\kappa)}\rho^{(2\kappa-2\mu)/(2-\kappa)
}\big)\| u\|_{E_{p,q,\beta}}.
$$
For $\epsilon\leq \rho^{2-\mu}$ this yields (here we use that $ \kappa\leq \mu<2$)
$$
 I_{\rho}\leq N\epsilon \|\partial_{t}u,
D^{2}u\|
_{ E _{p,q,\beta}}+N\epsilon^{-\mu/(2-\mu)}\| u\|_{E_{p,q,\beta}}.
$$

In the remaining case $\rho^{2-\mu}<\epsilon$. In that case for $\zeta\in
C^{\infty}_{0}((-1,1)\times B_{2})$   such that $\zeta
=1$ on $C_{1}$ we have by Lemma \ref{lemma 6,18.1} that
\begin{equation}
                     \label{8,27.2}
I_{\rho}\leq \rho^{2-\mu} \rho^{\beta+\gamma-2}
\dashnorm D(\zeta u)(0,\cdot)\|_{L_{r}(B_{\rho})}\leq N\epsilon
\|\partial_{t}(\zeta u),D^{2}(\zeta u)
\|_{\dot E_{p,q,\beta}}.
\end{equation}
Owing to $\beta\leq d/p+2/q$,
 the last norm here is easily shown to be less than
$$
N\|\partial_{t} u,D^{2}u
\|_{  E_{p,q,\beta} }
+N\| u\|_{  E_{p,q,\beta}}.
$$ 
Therefore, \eqref{6,20.1} holds in this case as well and this proves
\eqref{6,6.1}. Estimate \eqref{8,27.1}
follows from \eqref{8,27.2} with
$\mu=2$ and $\epsilon=1$.
The theorem is proved. \qed

\mysection {Solvability of parabolic equations and weak uniqueness of solutions of \protect\eqref{11.29.20}}
                 \label{section 3.27.3}

Assume \eqref{3.21.01}.
Fix some $\rho_{a},\rho_{b} \in(0,1]$.
The (small) parameters $\check a=\check a(d,\delta,p,q)>0$ and 
$$
\check b  :=\check b  (d,\delta, p, 
q, \rho_{a}, \beta  ) \in(0,1] 
$$  
below are taken from Theorem 3.5  of \cite{2} with $\alpha=1/2$.
\begin{assumption}
          \label{assumption 12.12.3}
 We have   
\begin{equation}
                                 \label{6.3.1}
 a^{\sharp}_{ \rho_{a}} =\sup_{\substack{\rho\leq\rho_{a}\\C\in\bC_{\rho}}}\dashint_{C}
\dashint_{C}|a ( x)-a(y)|
\,dx dy \leq \check a.
\end{equation}

\end{assumption}

\begin{assumption}
                      \label{assumption 3.14.1}
The function $c$ is bounded and 
there exists $p_{b}\in[p\beta,\infty)$ such that
\begin{equation}
                           \label{3.14.2}
\hat b_{\rho_{b }} =\sup_{r\leq\rho_{b }}r
\sup_{B\in \bB_{r}} 
\dashnorm b \|_{L_{p_{b} }(B)}\leq   \check b .
\end{equation}
\end{assumption}

In Assumption \ref{assumption 3.14.1} 
we have
$b\in L_{p_{b},\loc}$. It might happen
the $b$ is not summable to a much higher
degree.
\begin{example}
                 \label{example 7,29.1}
Let $d=d'+d''\geq3$, $d'\geq2$, and $|b(x)|\leq
1/|x'|$, where $x'=(x^{1},...,x^{d'})$.
Then $b\in E_{p,1}$ for any $p\in[1,d')$,
but $b\not\in L_{d',\loc}$.

\end{example}

Recall that according to Theorem \ref{theorem 8,23.3} the operator $\cL$
is bounded as an operator from $E^{1,2}_{p,q,\beta}$ into $E_{p,q,\beta}$.

Here is the main result,
Theorem 3.5,  of \cite{2} adjusted
to our needs in which $\nu
=\nu(d,\beta,p,q,\rho_{b})>0$ is defined in Remark 2.2 of \cite{2}.  

\begin{theorem}
                       \label{theorem 10.27.1}
 Let $c$ be a bounded Borel function
on $\bR^{d+1}$.
Under the above assumptions there exists 
$$
\check \lambda_{0}=\check \lambda_{0}
(d,\delta, p, 
q, \rho_{a}, \beta ,     
\sup|c|)\geq1
$$
 such that,  for  any
$\lambda\geq  \rho_{b}^{-2}\check \lambda_{0} $   and $f\in E_{p,q,\beta  }$, there exists
a unique $E^{1,2}_{p,q,\beta  }$-solution $u$
of 
\begin{equation}
                      \label{6,13.3}
\partial_{t}u+\cL u-(c+\lambda) u=-f.
\end{equation}
 Furthermore, there exists
a constant $N$ depending only on  $d$, $\delta$, $p$, $q$, $\rho_{a}$,  $\beta $, 
   $ 
\sup|c|$,    
such that 
\begin{equation}
                                \label{10.28.01}
\|\partial_{t} u,D^{2}u, \sqrt\lambda Du, \lambda u\|_{E_{p,q,\beta}}
\leq   N \nu^{-1}\|f\|_{E_{p,q,\beta}}.    
\end{equation}
In addition, if $f(t,x)=0$ for all $x$
and $t\geq T$ for some $T$, then
$u(t,x)=0$ for $t\geq T$ and all $x$.
\end{theorem}

\begin{remark}
              \label{remark 6,13.1}
Here we discuss what happens if
$f$ belongs to a few different 
$E_{p,q,\beta}$.
Let $(p_{i},q_{i},\beta_{i})$,
$i=1,2$, satisfy \eqref{3.21.01}.
Since $E_{p_{i},q_{i},\beta_{i}  }\subset
E_{p_{0},q_{0},\beta_{0}}$,
where $(p_{0},q_{0},\beta_{0})
=(p_{1}\wedge p_{2},q_{1}\wedge q_{2},\beta_{1}\vee \beta_{2})$, and
$(p_{0},q_{0},\beta_{0})$ {\em satisfies\/}
\eqref{3.21.01}, the solutions
of \eqref{6,13.3} in $E^{1,2}_{p_{i},
q_{i},\beta_{i}}$ are independent of $i$ as long as  for  
$(p_{i},q_{i},\beta_{i})$, $i=0,1,2$,
conditions \eqref{6.3.1} and \eqref{3.14.2}    hold with the corresponding $ \theta,\check b$
and $\lambda\geq\rho_{b}^{-1}(
\check\lambda^{0}_{0}\vee
\check\lambda^{1}_{0}\vee
\check\lambda^{2}_{0})$ with  $\check\lambda^{i}_{0}$ corresponding to 
$(p_{i},q_{i},\beta_{i})$
and $f\in \cap_{i} E_{p_{i},q_{i},\beta_{i}}$.

\end{remark}

\begin{remark}
                     \label{remark 10.27.2}   
The unique solution $u$ from Theorem \ref{theorem 10.27.1}
possesses the following properties

a) obviously, $u\in W^{1,2}_{p,q,\loc}$;

b)   by Lemma \ref{lemma 3.27.5}, we have
$Du\in L_{r,s,\loc} $, where $(r,s)=  (\beta -1)^{-1}\beta (p ,q )$;

c) by Lemma  
\ref{lemma 3.27.40}   $u$ is bounded and   continuous if $\beta<2$.

\end{remark}

\begin{remark}
                   \label{remark 5,28.1}

Theorem \ref{theorem 10.27.1} is applicable
in the case of $ c,f$   independent of $t$ ($a,b$ are independent of $t$
by assumption),
and owing to uniqueness yields solutions that
are also independent of $t$. This leads to
the ``elliptic'' counterpart, (recall that
$E_{p,\beta }$ and $E^{ 2}_{p,\beta  }$
are the subsets of $E_{p,q,\beta }$ and $E^{1,2}_{p,q,\beta }$ consisting of time-independent functions) which says that  Theorem
\ref{theorem 10.27.1} holds true if we replace
$E_{p,q,\beta }$ and $E^{1,2}_{p,q,\beta }$ 
with $E_{p,\beta }$ and $E^{ 2}_{p,\beta }$, respectively. 
\end{remark}

In the following lemma as in similar
situations elsewhere by $u(x)$
we mean the value at $x$ 
of the continuous modification
of $u\in E^{2}_{p,\beta}\subset E^{1,2}_{p,q,\beta}$,
which exists by Lemma \ref{lemma 3.27.40}.

\begin{lemma}
                 \label{lemma 7,5.2}
Let $\beta<2$ and $x_{t}$ be an $E_{p, \beta}$-admissible solution   of 
\eqref{11.29.20} with $x=0$. Let $\lambda>0$, $f$
be Borel bounded functions on $\bR^{d}$,   and 
let $u\in E^{2}_{p,\beta}$ be
a solution of $\cL u-\lambda u=-f$.
Then for any $t\geq0$
\begin{equation}
                      \label{7,6.3}
u(x_{t})e^{-\lambda t}=E\Big(\int_{0}^{\infty}e^{-\lambda (t+s)}f(x_{t+s})\,ds \mid \cF_{t}\Big).
\end{equation}
\end{lemma}

Proof. According to Theorem \ref{theorem 3.16.1}, It\^o's formula is applicable
to $ u( x_{t})e^{-\lambda t} $, which yields   that
on the set $\{\tau_{R}>t\}$ 
for $T>t$ we have
$$
u(x_{t})e^{-\lambda t}=
E\Big(u(x_{ T\wedge \tau_{R}})
e^{-\lambda (T\wedge \tau_{R})}
\mid \cF_{t}\Big)
$$
$$
+E\Big(\int_{t}^{T\wedge \tau_{R}}
e^{-\lambda s}
 f(x_{ s})\,ds \mid \cF_{t}\Big).
$$
By sending $R,T\to \infty$ 
and, taking into account
that $\tau_{R}\to\infty$, $ f( x_{t})$ is bounded and $\lambda>0$,
we come to
\eqref{7,6.3} and the lemma is proved. \qed

{\bf Proof of Theorem \ref{theorem 7,6.1}}. 
Let $x_{t}$ be $E_{p,\beta}$-admissible.
Then, for continuous $f$ with compact support and $\lambda$
large enough there exists a unique $u\in E^{2}_{p,\beta}$ satisfying $\lambda u-\cL u=f$. It follows from Lemma \ref{lemma 7,5.2} that, for any $n=0,1,...$, 
  continuous $f_{i}$, $i=0,...,n$, with compact support and $0=t_{0}<...<t_{n}<t$,  we have
$$
 E\prod_{i=0}^{n}f_{i}(x_{t_{i}}) u(x_{ t_{n}})
= E\prod_{i=0}^{n}f_{i}(x_{t_{i}})\int_{0}^{\infty}e^{-\lambda s}
f(x_{t_{n}}+s)\,ds
$$
\begin{equation}
                   \label{7,6.4}
= \int_{0}^{\infty}e^{-\lambda s}E\prod_{i=0}^{n}f_{i}(x_{t_{i}})
f(x_{t_{n}}+s)\,ds
\end{equation}

On the right we have the Laplace
transform of a function continuous
on $(0,\infty)$. Recall that
by knowing the Laplace transform
for all $\lambda$ large allows one to reconstruct uniquely the underlying function almost everywhere or everywhere where it is continuous.
For $n=0$  
\eqref{7,6.4} shows that $Ef(x_{s})$
is uniquely defined by $f,s,a,b$.
Assuming that
$$
E\prod_{i=0}^{n}f_{i}(x_{t_{i}}) 
$$
is uniquely defined by $f_{i},t_{i},a,b$ and using \eqref{7,6.4}
we conclude that the same is true with $n+1$ in place of $n$. This
implies that the finite dimensional
distribution of any $E_{p,\beta}$-admissible solution of
\eqref{11.29.20} with $x=0$ is the same
as well and the theorem is proved.
\qed
 
{\bf Conflict of interest}. The author have
no conflict of interest. There is no data in the paper.

\end{document}